%% file: main.tex
    \newcommand{\forArxiv}[1]{#1}  
    \newcommand{\forJournal}[1]{}      
\forJournal{
	\documentclass[final]{l4dc2024}
}
\forArxiv{
	\documentclass{article}
}

\def\texpath{tex}
\def\figpath{figures}
\input{\texpath/packages.tex}
\forJournal{
	\title[Minimax dual control]{Minimax dual control with finite-dimensional information state}
	\author{%
	 \Name{Olle Kjellqvist} \Email{olle.kjellqvist@control.lth.se}\\
	 \addr Department of Automatic Control, Lund University, Sweden
	}
}
\forArxiv{
	\title{Minimax dual control with finite-dimensional information state}
	\author{Olle Kjellqvist\thanks{Department of Automatic Control, Lund University, Sweden}}
}
\input{\texpath/definitions.tex}

\begin{document}
\maketitle

\begin{abstract}%
	This article considers output-feedback control of systems where the function mapping states to measurements has a set-valued inverse. We show that if the set has a bounded number of elements, then minimax dual control of such systems admits finite-dimensional information states. We specialize our results to a discrete-time integrator with magnitude measurements and derive a surprisingly simple sub-optimal control policy that ensures finite gain of the closed loop. The sub-optimal policy is a proportional controller where the magnitude of the gain is computed offline, but the sign is learned, forgotten, and relearned online.

The discrete-time integrator with magnitude measurements captures real-world applications such as antenna alignment, and despite its simplicity, it defies established control-design methods. For example, whether a stabilizing linear time-invariant controller exists for this system is unknown, and we conjecture that none exists.
\end{abstract}
\forJournal{
}

\section{Introduction}
\label{sec:introduction}
\input{\texpath/introduction.tex}

\section{Minimax dual control}
\label{sec:minimax}
\input{\texpath/minimax.tex}
\section{Magnitude Control}
\label{sec:magnitude}
\input{\texpath/magnitude.tex}
\section{Conclusion}
\label{sec:conclusion}
\input{\texpath/conclusion.tex}
\forJournal{
	\acks{I was introduced to the problem of output feedback from magnitude measurements by Bo Bernhardsson at the start of my Ph.D. studies. This manuscript has benefited from lengthy discussions with the author's colleagues at the Department of Automatic Control. In particular, I want to thank Anders Rantzer and Tore Hägglund for their encouragement and feedback along the way, Venkatraman Renganathan for reviewing an earler version of this manuscript and Jonas Hansson for tackling my conjecture with much enthusiasm.

	This project has received funding from the European Research Council (ERC) under the European Union's Horizon 2020 research and innovation programme under grant agreement No 834142 (ScalableControl).
	The author is a member of the ELLIIT Strategic Research Area at Lund University. 
	}
}
\forArxiv{
	\section*{Acknowledgements}
	I was introduced to the problem of output feedback from magnitude measurements by Bo Bernhardsson at the start of my Ph.D. studies. This manuscript has benefited from lengthy discussions with the author's colleagues at the Department of Automatic Control. In particular, I want to thank Anders Rantzer and Tore Hägglund for their encouragement and feedback along the way, Venkatraman Renganathan for reviewing an earler version of this manuscript and Jonas Hansson for tackling my conjecture with much enthusiasm.

	This project has received funding from the European Research Council (ERC) under the European Union's Horizon 2020 research and innovation programme under grant agreement No 834142 (ScalableControl).
	The author is a member of the ELLIIT Strategic Research Area at Lund University. 
}
\bibliography{references}

\end{document}

%% file: tex/packages.tex
\forArxiv{
	\usepackage{amsmath}
	\usepackage{amssymb}
	\usepackage{amsthm}
	\usepackage{natbib}
	\usepackage{graphicx}
	\usepackage{url}
	\usepackage{algorithm2e}
	\usepackage{subfigure}
	\usepackage[margin=1.2in]{geometry}
	\usepackage{natbib}
	\bibliographystyle{plainnat}
}
\usepackage{times}
\usepackage{xcolor}
\usepackage{booktabs}
\usepackage{verbatim}
\usepackage{enumitem}

%% file: tex/definitions.tex
\forArxiv{
\newtheorem{theorem}{Theorem}
\newtheorem{lemma}[theorem]{Lemma}

\newtheorem{proposition}[theorem]{Proposition}
\newtheorem{remark}[theorem]{Remark}
\newtheorem{example}{Example}
}

\newcommand{\R}{\mathbb{R}}

\newcommand{\Hinf}{\mathcal{H}_\infty}

\DeclareMathOperator{\sign}{sign}
\DeclareMathOperator{\B}{\mathcal{B}} 

\newcommand{\seq}[3]{#1_{#2:#3}}

\newcommand\defeq{\triangleq}

\newcommand{\W}{\mathcal{W}}
\newcommand{\Y}{\mathcal{Y}}
\newcommand{\X}{\mathcal{X}}
\newcommand{\U}{\mathcal{U}}

\newtheorem{assumption}{Assumption}
\newtheorem*{excont}{Example \continuation}
\newcommand{\continuation}{??}
\newenvironment{continueexample}[1]
 {\renewcommand{\continuation}{\ref{#1}}\excont[continued]}
 {\endexcont}

\newcommand{\br}{\mathbf{r}}

\newcommand{\one}{\mathbf{1}}

\newcommand{\bmat}[1]{\begin{bmatrix} #1 \end{bmatrix}}
\newcommand{\tran}{{\mkern-1.5mu\mathsf{T}}}

%% file: tex/introduction.tex
\begin{figure}
		\centering
		\subfigure{%
		\label{fig:receiver}%
		\includegraphics[width=.3\linewidth]{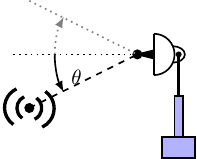}
	}
		\qquad
		\subfigure{%
		\label{fig:pendulum}%
		\includegraphics[width=.4\linewidth]{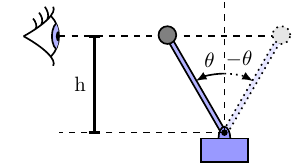}
	}
	\caption{%
		Examples of physical systems where the sign of the state is ambiguous:
		The left figure illustrates a receiver with an uncertain and potentially non-stationary source location. 
		The objective is to adjust the receiver's position to an angle that maximizes signal intensity. 
		Typically, the receiver's radiant sensitivity is symmetric relative to deviations from the incidence angle. 
		The right figure shows an inverted pendulum, which is regulated by monitoring the pendulum's height.
}
\end{figure}

This article concerns ouput feedback control of discrete-time systems whose measurement equations have a bounded number of solutions.
As a prototype example, we consider the discrete-time integrator, where the controller only has access to the magnitude of the state.
The state $x_t$, the control signal $u_t$, and disturbance $w_t$ are real-valued scalars.
The system is described by the recursion 
\begin{equation}\label{eq:state}
	x_{t+1} = x_t + u_t + w_t.
\end{equation}
We consider causal control policies, $\mu$, that map measurements of the state magnitude
\begin{equation}\label{eq:measurement}
	y_t = |x_t|
\end{equation}
to control signals
\begin{equation}\label{eq:admissible}
	u_t = \mu_t(y_0, y_1, \ldots, y_t, u_0, \ldots, u_{t-1}).
\end{equation}
The uncertain sign in \eqref{eq:measurement} captures some of the difficulties that may arise when optimizing a system based on measurements of some (locally) convex or concave performance quantity, as in Figure~\ref{fig:receiver}.
The problem is also closely related to stabilizing an inverted pendulum by feedback from height measurements rather than angular measurements, as in Figure~\ref{fig:pendulum}.
This plain-looking problem captures a surprising amount of complexity:
\begin{enumerate}
	\item \emph{Exploration vs. exploitation.}
		The more effectively we control the system, the less confident we become about the state's sign. If the system ever reaches $y = 0$, the state's sign information is lost.
	\item \emph{No stabilizing linear time-invariant controller.}
		Previous work report no stabilizing linear time-invariant controller for the system~\eqref{eq:state}--\eqref{eq:admissible}~\cite{Rosdahl2020Dual, Alspach1972Dual} and the system cannot be stabilized by proportional feedback\footnote{%
			A linear time-varying controller can stabilize the system.
			For example, $u_t = (-1)^t y_t$ will ensure $x_t = 0$ for all $t \geq 2$, for any $x_0$ and $w_t = 0$.
		}.
		This author \emph{conjectures} that there exists no finite-dimensional linear time-invariant controller that stabilizes the system.
	\item \emph{Extended Kalman filter.} 
		The extended Kalman filter (EKF) is a popular algorithm for estimating a nonlinear system's state, often coupled with certainty-equivalence control.
		However, the measurement equation~\eqref{eq:measurement} is not differentiable at $x=0$, and the EKF is not directly applicable.
		One may substitute the measurement equation with $y_t = x_t^2$ to recover differentiability, but this substitution results in an unobservable linearization.
	\item \emph{Myopic Controller.}
		The Myopic controller~\cite{Wittenmark1995Dual} associated with minimizing the current cost $x_t^2 + u_t^2$ is not stabilizing.
\end{enumerate}

In this article, we will design a control policy~\eqref{eq:admissible} that ensures that the induced $\ell_2$-gain from $w$ to $(x, u)$ is less than some positive quantity $\gamma$.
That is, the inequality
\begin{equation}\label{eq:finite_gain}
	\sum_{t=0}^N (x_t^2 + u_t^2) \leq \gamma^2 \sum_{t=0}^N w_t^2 + \beta(x_0)
\end{equation}
must be fulfilled for all $N \geq 0$, real-valued function $\beta$ and realizations $\seq{w}{0}{N} := w_0, w_1, \ldots, w_N$ of the disturbance sequence.
The condition~\eqref{eq:finite_gain} generalizes the classical $\Hinf$-norm for linear systems. 
The function $\beta$ is called a bias term and is used to capture the effect of the initial state.
The \emph{small-gain theorem} provides sufficient conditions for robust stability against feedback perturbations with induced $\ell_2$-norm less than $\gamma^{-1}$.
We refer the reader to~\cite[Chapter 5]{Khalil2002Nonlinear} for a detailed discussion on finite-gain stability and the small-gain theorem.
Surprisingly, we will see that it is possible to compress the observed output trajectory $(\seq{y}{0}{t}, \seq{u}{0}{t-1})$ using two recursively computed quantities $r^+_t$ and $r^-_t$.
These quantities correspond to the smallest feasible disturbance trajectory compatible with the observed outputs and $\sign(x_t) = 1$ or $\sign(x_t) = -1$.
Together with $y_t$, they make a \emph{sufficient statistic} for optimal control of a corresponding dynamic game.

The quantities follow the recursions
\begin{equation}\label{eq:rs}
	\begin{aligned}
		r^+_{t+1} & = y_t^2 + u_t^2 - \gamma^2\max\{r^+_t + (y_{t+1} - u_t - y_t)^2, r^-_t + (y_{t+1} - u_t + y_t)^2\}, \\
		r^-_{t+1} & = y_t^2 + u_t^2 - \gamma^2\max\{r^+_t + (-y_{t+1} - u_t - y_t)^2, r^-_t + (-y_{t+1} - u_t + y_t)^2\}.
\end{aligned}
\end{equation}
In Section~\ref{sec:magnitude}, we will show that these quantities are sufficient for ensuring bounded $\ell_2$ gain and summarize our conclusions about the magnitude control problem in Proposition~\ref{prop:magnitude}.
\begin{proposition}\label{prop:magnitude}
An admissible policy $\mu$ exists that ensures $\ell_2$-gain smaller than $\gamma$ if, and only if, it is achievable with a policy of the form $u_t = \eta_t(y_t, r^+_t, r^-_t)$.
Further, the controller $\eta(y, r^+, r^-) = 0.7\sign(r^- - r^+)y$ achieves $\ell_2$-gain less than $4$.
\end{proposition}
We remark that $\eta$ is admissible as $r^+_t$ and $r^-_t$ are functions of previous measurements and control signals.
Via substitution, one can recover $\mu$.

\subsection{Related work}
\paragraph{Adaptive control}
From the adaptive control perspective, system \eqref{eq:state}, \eqref{eq:measurement} could be interpreted as a linear system with uncertain time-varying parameters.
Several methods are described in excellent textbooks like~\cite[Chapter 6.7]{Goodwin2009Adaptive} that apply uncertain linear time-varying systems.
However, these methods rely on a separation of time scales between the state dynamics, the parameter adaptation, and the parameter variation.
Hence, we can not expect these methods to work well in our case~\cite{Anderson2008Challenges}. 
Nonlinear stochastic control theory provides a framework that can, in principle, handle fast parameter variation and large uncertainties, and our problem fits well with the methodology of \emph{dual control}~\cite[Chapter 7]{Wittenmark1995Dual}. 

Stochastic dual control has been applied to various problems with uncertain gain, as demonstrated in \cite{Astrom1986Dual, Dumont1988Wood, Allison1995Dual}.
\cite{Alspach1972Dual} considered control of an integrator based on noisy measurements of the square of the magnitude.
The noise was assumed Gaussian, and the author proposed approximating the information state by a sum of Gaussians.
\cite{Rosdahl2020Dual} considered a noisy version of the problem in this article but from a stochastic dual control perspective.
The authors proposed to approximate the information state by a neural network.

\paragraph{Learning-to-control}
Lately, there has been a surge of interest in learning to control linear systems.
Much of the work concerns the sample complexity of learning optimal controllers of linear time-invariant systems.
For example, \cite{Dean2018Regret, Mania2019Certainty} concerns quadratic performance objectives and additive stochastic noise, \cite{Chen2021Blackbox} adapts the theory of online convex optimization~\cite{Hazan2023Introduction} to unknown linear time-invariant systems with bounded disturbances.
\cite{Yu2023Online} proposed a method to control slowly varying linear systems with unknown parameters belonging to a polytope perturbed by bounded disturbances using convex body chasing.

\paragraph{Minimax control}
Minimax control for uncertain systems was introduced in the Ph.D. thesis of~\cite{Witsenhausen1966Minimax}. Information states, or sufficient statistics, for optimal control for output feedback minimax control, was discussed in \cite{Bertsekas1973Sufficient} based on Bertsekas's Ph.D. thesis. 
The game-theoretic formulation of $\Hinf$-control~\cite{Basar2008Hinf} is a special case of minimax control, and the information state formulation was derived for nonlinear systems in \cite{James1995Robust} demonstrating that, in general, the information state is infinite-dimensional.
The term minimax adaptive control was introduced in~\cite{Basar1994Minimax}. 
Recently,~\cite{Rantzer2021L4DC} proposed a minimax adaptive controller for uncertain linear systems with perfect state measurements.
The uncertainty was assumed to belong to a finite, known set.
The author proposed a finite-dimensional information state related to the empirical covariance matrix of the current state, previous state, and previous control signal.
This author extended Rantzer's results to scalar linear systems with noisy measurements in~\cite{Kjellqvist22Learning}.
Recently, \cite{Renganathan2023Online} studied the regret of Rantzer's controller for linear systems with energy-bounded disturbances.

\subsection{Contributions}
This article identifies a class of systems where the minimax dual controller admits a finite-dimensional information state.
The information state admits recursive computation, and Theorem~\ref{thm:verification} shows the equivalence between the minimax dual control problem and an information-state dynamic programming problem.
We also provide a dissipativity interpretation in Theorem~\ref{thm:approximation}.
The proofs of Theorems~\ref{thm:verification} and~\ref{thm:approximation} are available in the ArXiv version of this article~\cite{Kjellqvist2023Minimax}\forJournal{ and omitted here due to space constraints}.
These results generalize Theorem 1 in~\cite{Rantzer2021L4DC} to a larger system class and specialize the results in~\cite{James1995Robust} to classes of systems where the information state iteration becomes explicit.
The explicit iteration results from the bounded number of solutions to the measurement equation~\eqref{eq:measurement} and can be exploited to obtain closed-form (suboptimal) solutions to the minimax dual control problem.
We specialize these results to the magnitude control problem in the introduction and prove Proposition~\ref{prop:magnitude} in Section~\ref{sec:magnitude}.

\subsection{Notation}
We use $\R$ to denote the set of real numbers, $\R^n$ means the set of $n$-dimensional real vectors, and $\R^{n \times m}$ means the set of $n \times m$ real matrices.
The vector of ones is denoted $\one$.
We use $\seq{y}{0}{N}$ as shorthand for the sequence $(y_0, y_1, \ldots, y_N)$.
For a matrix $A \in \R^{n \times m}$, we denote the transpose by $A^\tran$.
For sets $A \subseteq S$ and $B \subseteq T$, and a function $f : S \to T$, the image of $A$ is denoted $f(A)$ and the preimage of $B$ is denoted $f^{-1}(B)$; the Cartesian product is denoted $S \times T$ and the $n$-ary Cartesian power $S = \underbrace{S \times S \times \ldots \times S}_{n \text{times}}$ is denoted $S^n$.
For vectors $v, v' \in \R^n$, the inequality $v \leq v'$ is understood component-wise, and for functions $f, g: S \to T$, the inequality $f \leq g$ means that $f(s) \leq g(s)$ for all $s \in S$ where $\leq$ is the partial order on $T$.
Strict inequalities are defined analogously.

%% file: tex/minimax.tex
This section introduces the minimax dual control problem, the information state, and dynamic programming.
By information state, we mean an auxiliary state variable that is computable by the controller, has a recursive expression in observed quantities and is sufficient to compute the optimal control policy and the associated cost.
For example, in the linear-quadratic Gaussian control problem, the information state is the conditional mean and covariance of the state given the observations---the Kalman filter estimate and the error covariance.
It is well known that the ``worst-case'' history is an information state for the minimax control problem, and dynamic programming with this information state is pretty well understood.
Unfortunately, this information state is generally infinite-dimensional and, therefore, impractical.
The main contribution of this section is to show that for our class of systems, the worst-case history admits a finite-dimensional representation.
This representation is, in itself, an information state.
We derive a verification and an approximation theorem for value iteration specific to this finite-dimensional representation.
\subsection{Problem formulation}

Let $f : \X \times \U \times \W \to \X$ and $h : \X \to \Y$ describe the dynamical system
\begin{equation}
	\begin{aligned}\label{eq:minimax_dynamics}
		x_{t+1} & = f(x_t, u_t, w_t) \\
		y_t & = h(x_t).
	\end{aligned}
\end{equation}
The control signal, $u_t \in \U$ is generated by a causal control policy $\mu_t : \Y^{t} \times \U^{t-1} \to \U$, where $\Y = h(\X)$ by
\begin{equation}\label{eq:JpiN}
	u_t = \mu_t(\seq{y}{0}{t}, \seq{u}{0}{t-1}).
\end{equation}
We call the tuple $\pi = (\mu_0, \mu_1, \ldots)$ a \emph{strategy} and the set of all such admissible strategies $\Pi$.
Consider the objective function as the ``worst-case'' sum of stage costs $l : \X \times \U \times \W \to \R$,
\begin{equation}\label{eq:J}
	J_\pi^N(y_0) \defeq \sup_{\seq{w}{0}{N}}\left\{\sum_{t=0}^{N} l(x_t,u_t, w_t) : \seq{w}{0}{N} \in \W^{N+1}, y_0 = h(x_0)\right\}.
\end{equation}
The goal of this section is to examine the minimax optimal control problem
\begin{equation}\label{eq:Jstar}
	J_\star(y_0) \defeq \inf_{\pi \in \Pi}\sup_N J_\pi^N(y_0).
\end{equation}

We make two crucial assumptions:
\begin{assumption}\label{ass:pos}
	For all $x \in \X, u \in \U$, $\sup_w l(x, u, w) \geq 0$.
\end{assumption}
The assumption that $\sup_w l(x, u, w) \geq 0$ implies monotonicity properties of $J_\pi^N$ in \eqref{eq:J} and, as we will see later, the value iteration.
\begin{assumption}\label{ass:finite}
	For any $y \in \Y$, the preimage $h^{-1}\{y\} \subset \X$ is an indexed set of at most $M$ elements.
\end{assumption}
This assumption relates to the dimensionality of the information state, or sufficient statistic, of the dynamic programming version of this problem.
Technically, the bound $M$ does not have to be known a priori, but we require the capability to enumerate all the solutions to $y_t = h(x_t)$ online.
At first glance, this assumption may appear overly limiting, but the following examples prove otherwise.

\begin{example}[Magnitude control of input-output models]\label{ex:magnitude}
	Consider controller design for the input-output system
	\begin{equation}\label{eq:IO}
		z_{t + 1} = -a_1 z_t - \cdots - a_d z_{t-d + 1} + b_1 u_t + \cdots + b_d u_{t-d + 1} + w_t,
	\end{equation}
	where the controller has access magnitude measurements $|z_0|, |z_1|, \ldots, |z_t|$ at time $t$.
	The system~\eqref{eq:IO} has a (nonminimal) state-space realization $x_{t+1} = Ax_t + Bu_t + Gw_t$, where
	\[
		x_t = \bmat{ z_t \\ \vdots \\ z_{t - d + 1} \\ u_{t - 1} \\ \vdots \\ u_{t - d + 1} }, \quad
		A = \bmat{ -a_1 & \cdots & -a_{d} & b_2 & \cdots & b_{d} \\
			1 & \\
			  & \ddots \\
			  & & 1 \\
			0 & \cdots & 0 & 0 & \cdots & 0 \\
			 & & & 1 \\
			 & & & & \ddots \\
			 & & & & & 1
		}, \quad
		B = \bmat{ b_1 \\ 0 \\ \vdots \\ 0 \\ 1 \\0 \\ \vdots \\ 0}, \quad
		G = \bmat {1 \\ 0 \\ \vdots \\ 0 \\ 0 \\ 0 \\ \vdots \\ 0}.
	\]
	Store the past $d-1$ inputs and outputs and define the augmented measurement 
	\[
		y_t = h(x_t) = (|z_t|, \ldots, |z_{t-n+1}|, u_{t-1}, \ldots, u_{t-n+1}).
	\]
	Then, the preimage
\[
	h^{-1}\{y_t\} = \{\pm |z_t|\} \times \cdots \times \{\pm |z_{t-d + 1}|\} \times \{u_{t-1}, \ldots, u_{t-d + 1}\}
\] 
has cardinality $2^{n}$, corresponding to the possible signs of the past measurements.
A first-order difference equation can model the integrator in the introduction, so $M = 2^1 = 2$, and the inverted pendulum (linearized around its equilibrium) by a second-order difference equation, for which $M = 2^2 = 4$.
\end{example}

\begin{example}[Linear system with uncertain dynamics]\label{ex:linear}
	Consider the linear system $x_{t+1} = Ax_t + Bu_t + w_t$ where $A, B$ are unknown matrices belonging to a finite set $\mathcal M$ of cardinality $M$.
	Then, the equivalent lifted system $x_t = (z_t, A_t, B_t)$ with
	\[
		\begin{aligned}
			A_{t+1} & = A_t, \quad B_{t+1}  = B_t, (A_0, B_0 ) \in \mathcal M\\
			z_{t+1} & = Az_t + Bu_t + w_t, \quad y_t = h(x_t) = z_t
		\end{aligned}
	\]
	satisfies Assumption~\ref{ass:finite} as $h^{-1}\{y_t\} = \{z_t\} \times \mathcal M$ has cardinality $M$.
\end{example}

\begin{example}[Finite state space]\label{ex:finite}
	If the state space $\X$ is finite, per definition $h^{-1}\{y\} \subseteq \X$ is finite.
\end{example}

\begin{remark}
	In our case $f$ and $h$ are given by \eqref{eq:state} and \eqref{eq:measurement} and the stage cost is $l(x_t, u_t, w_t) = x_t^2 + u_t^2 - \gamma^2 w_t^2$. 
	The states, observations and inputs take values in $\X = \R, \U = \R, \Y = \R_{\geq 0}, \W = \R$.
	The finite-gain condition~\eqref{eq:finite_gain} then correspond to $J_\star(y_0)$ being bounded.
	If not for the nonlinearity $h(x_t) = |x_t|$, it would be equivalent to the standard dynamic game formulation of $\Hinf$ suboptimal control~\cite{Basar2008Hinf}, rather it can be seen as a special case of nonlinear $\Hinf$ output feedback control~\cite{James1995Robust}.
\end{remark}

\subsection{An information state}
Following previous work \cite{Witsenhausen1966Minimax, Bertsekas1973Sufficient, James1995Robust, Basar2008Hinf} we consider the ``worst-case history'', $\rho_t$, that is compatible with the observations $\seq{y}{0}{t-1}$ and inputs $\seq{u}{0}{t-1}$ up to time $t-1$ reaching the state $x$ at time $t$:
\begin{multline}\label{eq:rho}
	\rho_t(x, \seq{y}{0}{t-1}, \seq{u}{0}{t-1})  \\
	\defeq \sup_{\seq{w}{0}{t-1} \in \W^{t-1}}\sup_{x_0 \in \X}\left\{ \sum_{\tau=0}^{t-1} l(x_\tau, u_\tau, w_\tau) : x_t = x, x_{\tau + 1} = f(x_\tau, u_\tau, w_\tau), y_\tau = h(x_\tau) \right\}.
\end{multline}
\begin{remark}
	We follow the convention that the supremum over the empty set is $-\infty$.
\end{remark}

The worst-case performance of a policy $\pi \in \Pi$, $J_\pi^N(y_0)$ can be expressed in terms of $\rho_t$ as
\begin{equation}\label{eq:J_rho}
	J_\pi^N(y_0) = \sup_{\seq{y}{0}{N}, x} \rho_{N+1}(x, \seq{y}{0}{N}, \seq{u}{0}{N}),
\end{equation}
where $\seq{u}{0}{N}$ is generated by $\pi$ and $\seq{y}{0}{N}$.
The functions $\rho$ are causal functions of the measurements and control signals and obey the forward dynamic programming, \cite{Magill1965Optimal}, recursion:
\begin{equation}\label{eq:rho_recursion}
	\rho_{t + 1}(x, \seq{y}{0}{t}, \seq{u}{0}{t}) = \sup_{\xi,w} \left\{ l(x, u_t, w) + \rho_t(\xi, \seq{y}{0}{t-1}, \seq{u}{0}{t-1}) : x = f(\xi_t, u_t, w), y_t = h(\xi) \right\}.
\end{equation}
Each step~\eqref{eq:rho_recursion} involves extremizing over the previous state and the disturbance dependent on the current state $x_t$, and in general, the computational complexity of evaluating $\rho_t$ grows with $t$.
However, for systems satisfying Assumption~\ref{ass:finite}, the set of feasible past states involved in \eqref{eq:rho_recursion} is restricted by the measurement trajectory.
To exploit this restriction, we split the computation of \eqref{eq:rho_recursion} into two steps: a correction step incorporating the observation $y_t$ and a prediction step after selecting $u_t$:
\begin{subequations}\label{eq:rho_computation}
	\begin{align}
		(\xi_t^i)_{i = 1}^M & =  h^{-1}\{y_t\} \label{eq:rho_computation:xi} \\
		r_t^i & = \max_{j} \rho_t(\xi_t^i, \seq{y}{0}{t-1}, \seq{u}{0}{t-1}) \label{eq:rho_computation:r}\\
		\rho_{t+1}(x, \seq{y}{0}{t}, \seq{u}{0}{t}) & = \sup_{i, w \in \W} \{l(\xi^i_t, u_t, w) + r_t^i : x = f(\xi_t^i, u_t, w) \}. \label{eq:rho_computation:rho}
	\end{align}
\end{subequations}

The intuition behind procedure~\eqref{eq:rho_computation} is that at time $t$, the realization of the state $x_t$ must belong to the $M$ solutions of $y_t = h(\xi)$. 
The value $r_t^i$ is the worst-case performance of the system up to time $t$ under the hypothesis that $x_t = \xi_t^i$ consistent with $\seq{y}{0}{t}$ and $\seq{u}{0}{t-1}$.
The prediction $\rho_{t+1}(x, \seq{y}{0}{t}, \seq{u}{0}{t})$ is the worst-case performance of the system up to time $t+1$ under the hypothesis that $x_{t+1} = x$ consistent with $\seq{y}{0}{t}$ and $\seq{u}{0}{t}$. 
The extremization~\eqref{eq:rho_computation:rho} includes two terms: the stage cost $l(x, u_t, w)$ capturing the cost of transition to $x_{t+1} = x$ from $x_t$ and the past performance $r_t^i$ under the hypothesis $x_t = \xi_t^i$.
The extremization is carried out over the hypotheses $\xi_t^1, \ldots, \xi_t^M$ and the disturbance $w$.
Define the update functions $g$ for the $M$-dimensional vector $\br_t = (r_t^1, \ldots, r_t^M)$ 
\begin{align}
	g_i(\br, \xi_+, \xi, u) & \defeq \sup_{j, w \in \W}\{ l(\xi^j, u, w) + r^j : \xi_+^i = f(\xi^j, u, w) \} \nonumber \\
	g(\br, y_+, y, u) & \defeq g_i(\br, h^{-1}\{y_+\}, h^{-1}\{y\}, u).\label{eq:g}
\end{align}
In the following proposition, we formalize the properties of the update functions $g$ and the sequence $\br$.
\begin{proposition}\label{prop:rprop}
Fix $N$, $i = 1, \ldots, M$, $y_0$, a strategy $\pi$ and let $\br_N$ be defined recursively by $\br_0 = 0$ and $\br_t = g(\br_{t-1}, y_t, y_{t-1}, u_{t-1})$.
Then
\[
	\rho_{t+1}(x, \seq{y}{0}{t}, \seq{u}{0}{t}) = \sup_{i,w} \left\{ l(x, u_t, w) + r_t^i : x = f(\xi_t^i, u_t, w) \right\}.
\]
Furthermore, \forJournal{1. There exists a sequence $\seq{w}{0}{N-1}$ such that $\max_i r_N^i \geq 0$. 
2. For fixed $y_t, y_{t-1} \in \Y$ and $u_t \in \U$, for $\br \leq \br'$ we have $g(\br, y_t, y_{t-1}, u_t) \leq g_{r}(\br', y_t, y_{t-1}, u_t)$.
3. $g(\br + \one c, y_t, y_{t-1}, u_t) = g(\br, y_t, y_{t-1}, u_t) + \one c$ for all $c \in \R$.
}%
\forArxiv{
\begin{enumerate}
	\item There exists a sequence $\seq{w}{0}{N-1}$ such that $\max_i r_N^i \geq 0$.
	\item For fixed $y_t, y_{t-1} \in \Y$ and $u_t \in \U$, for $\br \leq \br'$ we have $g(\br, y_t, y_{t-1}, u_t) \leq g_{r}(\br', y_t, y_{t-1}, u_t)$.
	\item $g(\br + \one c, y_t, y_{t-1}, u_t) = g(\br, y_t, y_{t-1}, u_t) + \one c$ for all $c \in \R$.
\end{enumerate}
}
\end{proposition}
\forArxiv{
\begin{proof}
	\emph{1.} follows directly from Assumption~\ref{ass:pos}. \emph{2.} follows from the monotonicity of the supremum operator. \emph{3.} follows from that for any function $f$ and set $Z$ $\sup_{z \in Z} \{f(z) + c\} = \sup_x\{f(z)\} + c$ for all $c \in \R$. 
	Finally, by recursion, the elements in $\br_t$ are equal to the ones in \eqref{eq:rho_computation:r}, thus for each $i = 1, \ldots, M$ equation~\eqref{eq:rho_computation:rho} holds with $(r_t^i)_{i = 1}^M = \br_t$.
\end{proof}
}

By Proposition~\ref{prop:rprop}, the worst-case history $\rho_N$, is sufficient to evaluate the objective $J_\pi^N(y_0)$. 
We will now study value iteration to minimize $\sup \rho_N$.
Consider the time evolutions of the measurements $y$ and representations $\br$:
\begin{subequations}\label{eq:information}
\begin{align}
	y_{t+1} & = v_t \\
	\br_{t+1} & = g(\br_t,v_t, y_t, u_t), \quad \br_0 = 0,  \label{eq:information:r_recursion}
\end{align}
\end{subequations}
where the next measurement, $v_t$, is considered an exogenous input.

The optimization problem~\eqref{eq:Jstar} can be expressed in terms of the worst-case history $\rho_N$ as
\begin{equation}\label{eq:information_state_problem}
	\inf_\eta \sup_{N, \seq{v}{0}{N-1}, x\in\X} \left\{ \rho_N(x, \seq{v}{0}{N-1}, \seq{u}{0}{N-1})\right\},
\end{equation}
where an information-state feedback policy generates $u_t$
\[
	u_t = \eta_t(\br_t, y_t).
\]
Define the set of information-state strategies $\widetilde \Pi$ as the set of strategies $\tilde \pi = (\eta_0, \eta_1, \ldots)$.
As $\br$ is a causal function of the measurements and control signals, so is $\eta_t$ (by composition) and $\widetilde \Pi \subset \Pi$.
	In other words, information-state feedback is admissible.
	The following examples illustrate the information-state recursions~\eqref{eq:information} for the systems in Examples~\ref{ex:magnitude} and~\ref{ex:linear}.
\begin{continueexample}{ex:magnitude}
	In this case, it is convenient to index the hypotheses $h^{-1}\{y_t\}$ by sequences of hypothetical signs, $s_{t}, \ldots, s_{t - d + 1}$ of the $d$ stored measurements $|z_t|, \ldots, |z_{t - d + 1}|$.
	The update simplifies significantly as the realizations of $z_t, \ldots, z_{t - d + 2}$ must remain unchanged between time steps $t$ and $t + 1$.
	Further, $w_t = z_{t + 1} + a_1 z_t + \ldots + a_dz_{t - d + 1} - b_1u_t + \ldots + b_du_{t-d + 1}$ is uniquelly determined by the state trajectory, so
	\begin{multline}
		r_{t + 1}^{s_{t + 1}, \ldots, s_{t - d + 2}} = \max_{s_{t - d + 1} = \pm 1}\Big \{
		l(s_t|z_t|, u_t, w) + r_t^{s_t, \ldots, s_{t - d + 1}}  \\
		: w = s_{t + 1}|z_{t + 1}| + a_1 s_t|z_t| + \ldots + a_d s_{t - d + 1}|z_{t - d + 1}| - b_1u_t + \ldots + b_du_{t-d + 1} \Big \}.
	\end{multline}
\end{continueexample}

\begin{continueexample}{ex:linear}
	Here, we index the hypotheses $h^{-1}\{y_t\}$ by the matrices $A_t, B_t$.
	The update becomes
	\[
		r_{t + 1}^{A_{t + 1}, B_{t + 1}} = \sup_{A_t, B_t}\left\{ l(x_t, u_t, x_{t+1} - A_t x_t - B_t u_t) + r_t^{A_t, B_t} : A_{t + 1} = A_t, B_{t + 1} = B_t \right\}.
	\]
	By assumption $(A_t, B_t) = (0, 0)$ for all $t$, so the update simplifies to $r_{t + 1}^{A, B} = l(x_t, u_t, x_{t+1} - A x_t - B u_t) + r_t^{A, B}$.
\end{continueexample}

\subsection{Value iteration}
Towards finding (sub)optimal solutions to~\eqref{eq:Jstar}, we introduce the Bellman operators $\B$ and $\B_u$ for functions $V : (\R\cup \{-\infty\})^M \times \Y \to \R$
\begin{equation}\label{eq:Bellman}
	\B V(\br, y) = \min_{u \in \U(y)}\overbrace{\max_{v \in \Y} \left\{ V(g(\br,v, y, u), v) \right\}}^{\B_u V(\br, y)}.
\end{equation}
and the value iteration
\begin{subequations}\label{eq:value_iteration}
	\begin{align}
		V_0(\br, y) & = \max_{i = 1, \ldots, M} \{ r^i\} \\
		V_{k + 1}(\br, y) & = \B V_k(\br, y).
	\end{align}
\end{subequations}

We are ready to state the main theoretical results, justifying the value iteration algorithm~\eqref{eq:value_iteration}.
\begin{theorem}\label{thm:verification}
	For the system~\eqref{eq:minimax_dynamics} under Assumptions~\ref{ass:pos} and \ref{ass:finite} and strategy class $\Pi$, the value~\eqref{eq:Jstar} is bounded for any $x_0 \in \X$ if, and only if, the sequence $V_0, V_1, \ldots$ defined in~\eqref{eq:value_iteration} is bounded.
	If bounded, the sequence converges to the optimal value function $V_\star$.
	The limit $V_\star$ is a fixed point of the Bellman operator~\eqref{eq:Bellman} and the value $J_\star(y_0) = V_\star(0, y_0)$.
	If the minimum in \eqref{eq:Bellman} is attained for some $u \in \U$ for all $y \in \Y$ and $\br$, then the policy $\eta_*(\br, y)$ defined as the minimizing argument in \eqref{eq:Bellman} satisfies $\B_{\eta_\star(\br, y)}V_\star(\br, y) = V_\star(\br, y)$ and the policy
	\[
		\mu_t(\seq{y}{0}{t}, \seq{u}{0}{t-1}) = \eta_\star(\br_t, y_t)
	\]
	is optimal for~\eqref{eq:Jstar}.
\end{theorem}
\forArxiv{
\begin{proof}
	For any fix $N \geq 0$, the quantity $\inf_{\pi \in \Pi} J_\pi^N(y_0)$ lower bounds $J_\star(y_0)$ due to Assumption~\ref{ass:pos}.
	By \eqref{eq:J_rho},
	\[
		\inf_{\pi \in \Pi}J_\pi^N(y_0) = \inf_{\pi \in \Pi}\sup_{x, w}\rho_{N+1}(x, \seq{y}{0}{N}, \seq{u}{0}{N}) 
		= \inf_{\pi \in \Pi} \sup_{i, \seq{v}{0}{N+1}}\{ r_{N + 1}^i : v_t \in \Y \}.
	\]
	By standard dynamic programming arguments, see for example~\cite[Chapter 1.6]{Bertsekas2005DP}, this is equal to
	\[
		\begin{aligned}
&	 \inf_{\seq{\mu}{0}{N-1}}\inf_{\mu_N} \sup_{\seq{v}{0}{N-1}} \sup_{v_N}\{V_0(\br_{N + 1}, y_{N+1}) : v_t \in \Y \} \\
&	= \inf_{u_0 \in \U}\sup_{v_0 \in \Y} \cdots \inf_{u_N \in \U} \sup_{v_N} V_0(\br_{N + 1}, y_{N+1})
	= V_{N+1}(0, y_0).
		\end{aligned}
	\]
	This proves that the sequence $V_0(0, y_0), V_1(0, y_0), \ldots$ is bounded if $J_\star(y_0)$ is bounded. 
	By assumption, this holds for all $y_0 \in \Y$.

	By induction, the value iteration is non-decreasing as $\B$ is monotone and $V_1 \geq V_0$ follows from assumption~\ref{ass:pos}.
	Further, $V_k(\br, y) \leq V_k(\max\{r_i\} \one, y) = V_k(0 , y) + \max r_i$, proving that $V_0, V_1, \ldots$ is bounded, and since it is monotone increasing, it converges to a limit $V_\star$.

	Assume that $V_0, V_1, \ldots$ is bounded towards proving the other direction. 
	Then $V_\star$ is well-defined and satisfies $V_\star \geq V_k$ for all $k$.
	Fix an arbitrary $\epsilon > 0$ and define a policy $\eta_t^\epsilon$ that chooses $u_t$ such that
	\[
		\B_{u_t}V_\star(\br_t, y_t) \leq V_\star(\br_t, y_t) + \frac{1}{2}\epsilon \left(\frac{1}{2}\right)^t.
	\]
	By the definition of the infimum, such a $u_t$ always exists.
	Then, by similar arguments as above, we have	
	\[
		J_\star(y_0) \leq \sup_N J_{\eta^\epsilon}^N(y_0) \leq \sup_N V_\star(0, y_0) + \frac{1}{2}\epsilon \sum_{t=0}^{N + 1} \left(\frac{1}{2}\right)^t = V_\star(0, y_0) + \epsilon.
	\]
	So we have $V_\star(0, y_0) \leq J_\star(y_0) \leq V_\star(0, y_0) + \epsilon$.
	As $\epsilon$ was arbitrary, we conclude $J_\star(y_0 ) = V_\star(0, y_0)$.
	If for any $y \in \Y$ and $\br$, the minimum in \eqref{eq:Bellman} is attained for some $u \in \U$, then we can pick $\epsilon = 0$ and conclude that the minimizing argument in \eqref{eq:Bellman} is optimal for \eqref{eq:Jstar}.
\end{proof}
}
\begin{theorem}[Approximation]\label{thm:approximation}
	For the system~\eqref{eq:minimax_dynamics} under Assumptions~\ref{ass:pos} and \ref{ass:finite} and strategy class $\Pi$,
	assume that there exists a function $\bar V : (\R \cup \{-\infty\})^M \times \Y \to \R$ and a strategy $\bar \pi = (\bar \eta, \bar \eta, \ldots) \in \widetilde \Pi$ such that $\bar V \geq V_0$ and
	\[
		\B_{\bar \eta(\br, y)} \bar V(\br, y) \leq \bar V(\br, y).
	\]
	Then the value iteration $V_0, V_1, \ldots$ is bounded, and $J_{\bar \mu}(y_0) \leq \bar V(0, y_0)$ for the policy
	
	\[
		\bar\mu_t(\seq{y}{0}{t}, \seq{u}{0}{t-1}) = \bar\eta(\br_t, y_t).
	\]
\end{theorem}
\forArxiv{
\begin{proof}
	By monotonicity of the Bellman operator, we have that $V_k \leq \bar V$ for all $k = 0, 1, \ldots$, implying that the value iteration $V_0, V_1, \ldots$ is bounded.
	Further,
	\[
		\begin{aligned}
			J_\star(y_0) & \leq \sup_N J_{\bar \pi}^N(y_0)
				     = \sup_N \sup_{\seq{v}{0}{N}}\{ V_0(\br_{N+1}, y_{N+1}) : v_t \in \Y \} \\
				     & \leq\sup_N \sup_{\seq{v}{0}{N}}\{ \bar V(\br_{N+1}, y_{N+1}) : v_t \in \Y \}
				      \leq \sup_N \bar V(\br_0, y_0) = \bar V(\br_0, y_0).
		\end{aligned}
	\]
\end{proof}
}

%% file: tex/magnitude.tex
We now apply the above results to the example in Section~\ref{sec:introduction}.
For any $y$, we denote $\xi^+ = y$ and $\xi^- = -y$.
Then $h^{-1}\{y\} = \{ \xi^+, \xi^- \}$.
We similarly index $\br = [r^+, r^-]$.
Then $g$ in \eqref{eq:g} becomes $
		g_s(\br, v, y, u) = y^2 + u^2 - \gamma^2\min\{r^+ + (sv -u -y)^2, r^- + (sv -u +y)^2\}
		$
for $s = \pm 1$.
\paragraph{Proof of Proposition~\ref{prop:magnitude}} 
By the above analysis, the quantities \eqref{eq:rs} correspond to \eqref{eq:information}, and the first statement in the proposition is a direct consequence of Theorem~\ref{thm:verification}.
Drawing inspiration from~\cite{Rantzer2021L4DC}, we parameterize an upper bound of the optimal value in the parameters $0 < p \leq q < \gamma^2$ by
\begin{equation}\label{eq:magnitude:barv}
	\bar V(\br, y) = \max \{py^2 + r^+, py^2 + r^-, q y^2 + (r^+ + r^-) / 2\},
\end{equation}
and a certainty equivalence policy
\begin{equation}\label{eq:etabar}
	\bar \eta(\br, y) = k\sign(r^- - r^+)y
\end{equation}

The following lemma relates the parameters of the value function approximation $p, q$ and $k$ to the $\ell_2$ gain of the closed loop.
\begin{lemma}\label{thm:magnitude}
	Given a quantity $\gamma > 0$, parameters $0 < p < q < \gamma^2$, $k \in \R$, and $\bar V$ as above.
	The certainty equivalency policy $\bar \eta$ in \eqref{eq:etabar} achieves an $\ell_2$-gain of at most $\gamma$ for the system \eqref{eq:state}--\eqref{eq:admissible} and an objective value smaller than $\bar V(0, |x_0|)$ for the decision problem~\eqref{eq:Jstar}, if
	\begin{equation}\label{eq:inequalities}
		\begin{aligned}
			p > 1 + k^2 + \frac{(1 - k)^2}{p^{-1} - \gamma^{-2}}, \quad 
			q > 1 + k^2 + \frac{(1 + k)^2}{p^{-1} - \gamma^{-2}}, \quad
			q > 1 + k^2 + \frac{1}{q^{-1} - \gamma^{-2}} - \gamma^2 k^2.
		\end{aligned}
	\end{equation}
\end{lemma}

\begin{figure}
		\centering
		\subfigure{%
		\label{fig:y}%
		\includegraphics{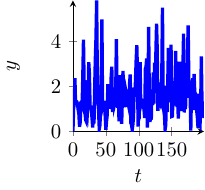}
	}
	\hfill
		\subfigure{%
		\label{fig:u}%
		\includegraphics{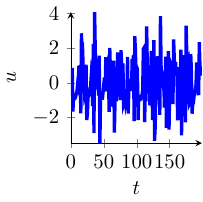}
	}
	\hfill
		\subfigure{%
		\label{fig:gain}%
		\includegraphics{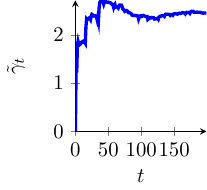}
	}
	\hfill
		\subfigure{%
		\label{fig:V}%
		\includegraphics{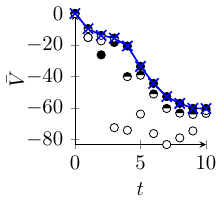}
	}

	\caption{%
		Figures~\ref{fig:y} and \ref{fig:u} contain plots of the outputs and control signal, respectively. 
		Figure~\ref{fig:gain} shows the empirical gain from $w$ to $(y, u)$ and Figure~\ref{fig:V} shows the value function approximation $\bar V = \max\{\bar V^{++}, \bar V^{+-}, \bar V^{--}\}$ defined in~\eqref{eq:magnitude:barv}.
		The black marks corresponds to $\bar V^{++}$, the half circles to $\bar V^{+-}$, the white marks to $\bar V^{--}$ and the blue crosses to $\bar V$.
		Note that the value function approximation is monotonically decreasing.
}
\label{fig:results}
\end{figure}

The values $\gamma = 4$, $p = 1.7$, $q = 7$ and $k = 0.7$ satisfy the conditions of Theorem~\ref{thm:magnitude} and a simulation with $w_t = sin(\pi t /10)$ is shown in Figure~\ref{fig:results}.
\forArxiv{
\begin{proof}
	Define
	\begin{equation}\label{eq:Vbars}
		\bar V^{++}(\br, y) \defeq py^2 + r^+, \quad \bar V^{--}(\br, y) \defeq py^2 + r^-, \quad \bar V^{+-}(\br, y) \defeq qy^2 + (r^+ + r^-) / 2.
	\end{equation}
	Then, for a fixed $u$, we have
	\[
		\begin{aligned}
			& \B_u \{\max\{\bar V^{++}(\br, y), \bar V^{--}(\br, y))\} 
			 = y^2 + u^2 + \max\left\{ 
			\frac{(y + u)^2}{p^{-1} - \gamma^{-2}} + r^+, \frac{(y - u)^2}{p^{-1} - \gamma^{-2}} + r^-
			\right\}.
		\end{aligned}
	\]
	Define for $i, j \in \{+, -\}$
	\[
		\alpha^{ij} \defeq \sup_{v \geq 0}\{qv^2 - \frac{\gamma^2}{2} \left( (v -u - iy)^2 + (-v - u - jy)^2 \right) + \frac{r^i + r^j}{2}.
	\]
	Then $\B_u \{\bar V^{+-}(\br, y)\} = y^2 + u^2 + \max_{i, j \in \{+, -\}} \alpha^{ij}$,	where, for $i \neq j$, 
	\[
		\alpha^{ii} = r^i - \gamma^2(u + iy)^2, \quad
		\max\{\alpha^{ij}, \alpha^{ji}\} = \frac{y^2}{q^{-1} - \gamma^{-2}} - \gamma^2 u^2 + \frac{r^i + r^j}{2}
	\]
	Let $l = \arg \max_{i \in \{+, -\}} \{r^i\}$, then by \eqref{eq:inequalities}, we have 
	\begin{multline*}
		\B_{-k_ly} \{\bar V(\br, y)\} = \max_{ij} \{\B_{-k_ly} \bar V^{ij}(\br, y)\} \leq \max\{\B_{-k_ly} \bar V^{il}(\br, y) -(r^l - r^j) / 2\} \\
		\leq \max\{ \bar V^{il}(\br, y)\} \leq \bar V(\br, y).
	\end{multline*}
	Therefore, by Theorem~\ref{thm:approximation}, the objective value is bounded from above by $\bar V(0, y_0)$.
\end{proof}
}

%% file: tex/conclusion.tex
This article demonstrated that output feedback minimax dual control possesses a finite-dimensional information state when the measurement equation has a finite number of solutions.
We applied this finding to the magnitude control of an integrator, resulting in a surprisingly simple sub-optimal control policy.
The controller is proportional, with the gain determined through hypothesis testing and updated online.
However, the results are limited to cases where the measurement equation has a finite number of solutions.
This restriction excludes scenarios where measurements are affected by real-valued sensor noise, which typically leads to an infinite-dimensional information state.

Future work will focus on extending these results to cases with noisy measurements, specifically where the dynamics are linear and uncertain but belong to a finite set.
Progress has already been made for scalar systems~\cite{Kjellqvist22Learning}, and the extension to multi-dimensional cases is currently under investigation.